\documentclass[aos,MSNbibl,nameyear,dvips]{arximspdf}

\doi{10.1214/12-AOS1069} 
\volume{40}
\issue{6}
\pubyear{2012}
\firstpage{3161}
\lastpage{3175}

\makeatletter

\newcommand{\rright}{\right}
\newcommand{\lleft}{\left}
\newtheorem{theorem}{Theorem}
\newtheorem{lemma}{Rule}
\newtheorem{corollary}{Corollary}
\newproclaim{example}{Example}

\makeatother

\begin{document}
\begin{frontmatter}

\title{A code arithmetic approach for quaternary code designs and its
application to $(1/64)${th}-fractions\thanksref{T1}}
\thankstext{T1}{Supported by National Science Council of Taiwan ROC
Grants 98-2118-M-001-028-MY2 and 100-2118-M-001-002-MY2.}
\runtitle{CA approach for QCD and applications to $(1/64)$th-fractions}

\begin{aug}
\author[A]{\fnms{Frederick K.~H.} \snm{Phoa}\corref{}\ead[label=e1]{fredphoa@stat.sinica.edu.tw}}
\runauthor{F.~K.~H. Phoa}
\affiliation{Academia Sinica}
\address[A]{Institute of Statistical Science\\
Academia Sinica\\
Taipei 11529\\
Taiwan\\
\printead{e1}} 
\end{aug}

\received{\smonth{1} \syear{2012}}
\revised{\smonth{11} \syear{2012}}

%
\begin{abstract}
The study of good nonregular fractional factorial designs has received
significant attention over the last two decades. Recent research
indicates that designs constructed from quaternary codes (QC) are very
promising in this regard. The present paper aims at exploring the
fundamental structure and developing a theory to characterize the
wordlengths and aliasing indexes for a general $(1/4)^p$th-fraction QC
design. Then the theory is applied to $(1/64)${th}-fraction QC
designs. Examples are given, indicating that there exist some QC
designs that have better design properties, and are thus more
cost-efficient, than the regular fractional factorial designs of the
same size. In addition, a~result about the periodic structure of
$(1/64)${th}-fraction QC designs regarding resolution is stated.
\end{abstract}

%


\end{frontmatter}

\section{Introduction}\label{sec1}

In many scientific researches and investigations, the interest lies in
the study of effects of many factors simultaneously. One may choose a
full factorial design which is able to estimate all possible level
combinations of factors, but it usually involves many unnecessary
trials. To be more cost-efficient, a fractional factorial design is
suggested. A good choice of fractional factorial design allows us to
study many factors with relatively small run size but enables us to
estimate a large number of effects.

Designs that can be constructed through defining relations among
factors are called regular designs, and all other designs that do not
possess this kind of defining relation are called nonregular designs.
\citet{WuHam00} and \citet{MukWu06} provide detailed
discussions on optimality criteria such as resolution and minimum
aberration for choosing fractional factorial\vadjust{\goodbreak} designs. Nonregular
designs have received particular attention in the past ten to twenty
years. The notions of resolution and aberration have been generalized
with statistical justifications to these designs; see \citet{DenTan99} and
\citet{TanDen99}. It is well recognized that although
nonregular designs have a complex aliasing structure, they can
outperform their regular counterparts with regard to resolution or
projectivity, and this is a major motivating force for the current
surge of interest in these designs. A comprehensive review on the
development of nonregular designs is referred to \citet{XuPhoWon09}.

A recent major development in nonregular two-level designs has been the
use of quaternary codes for their simple construction, and the
resulting two-level designs are generally called QC designs. \citet{XuWon07}
pioneered research on QC designs and reported theoretical
as well as computational results. \citet{PhoXu09} investigated the
properties of quarter-fraction QC designs. In addition to giving
theoretical results on the aliasing structure of such designs, they
constructed optimal quarter-fraction QC designs under several criteria.
\citet{Zhaetal11} introduced a trigonometric representation for the
study of QC designs and successfully derived the properties of
$(1/8)${th}- and $(1/16)${th}-fractions QC designs. The optimal
$(1/8)${th}- and $(1/16)${th}-fractions QC designs under maximum
resolution criterion were reported in \citet{PhoMukXu12}.

The present paper aims at exploring the fundamental structure and
developing the underlying theorems of a general QC design. In Section~\ref{sec2}
we recall some concepts about the design construction method via
quaternary codes. Then we introduce some new notation that is related
to wordlengths and aliasing indexes of words. This new notation
provides clear and simple presentations for theorems and examples in
the later sections. Section~\ref{sec3} contains some rules and corollaries about
the structure of QC designs. One can derive the wordlengths and
aliasing indexes of a word in a general QC design using these rules. In
addition, two theorems are stated about the structure of the $k$-equation
and their necessary and sufficient conditions. These theorems are
applied in Section~\ref{sec4}, leading to a theorem about the properties of
$(1/64)${th}-fraction QC designs. An example demonstrates the use of
the theorem to derive the generalized resolutions and generalized
wordlength patterns of QC designs. Based on the properties of the
derived classes of QC designs, the structure periodicity of
$(1/64)${th}-fraction QC designs with high resolution is suggested.
The proofs of these theorems are given in the last
section.\looseness=-1

\section{Definitions and notation}\label{sec2}

We recall some concepts in \citet{PhoXu09} here. A quaternary code
takes on values from $Z_4=\{0,1,2,3\}$. Let $G$ by an $n \times m$
generator matrix over $Z_4$. All possible linear combinations of the
rows in $G$ over $Z_4$ form a quaternary linear code, denoted by $C$.
Then each $Z_4$ entry of $C$ is transformed into two binary codes in
its binary image $D=\phi(C)$ via the Gray map, which is defined as follows:
\[
\phi\dvtx %
\matrix{ 0\to(1,1), & 1\to(1,-1), & 2\to(-1,-1), & 3
\to(-1,1). }\vadjust{\goodbreak} %
\]
Note that $D$ is a binary $2^{2n} \times2m$ matrix or a two-level
design with $2^{2n}$ runs and $2m$ factors.

In general, for highly-fractionated QC designs, we consider an $n
\times(n+p)$ generator matrix $G=(V,I_n)$, where $V=(\vec{v}_1,\ldots
,\vec{v}_p)$ is a matrix over $Z_4$ that consists of $p$ vectors of
lengths $n$ and $I_n$ is an $n \times n$ identity matrix. It leads to a
two-level design $D$ with $2^{2n}$ runs and $2n+2p$ factors, that is,
$D=(d_1,\ldots,d_{2p},d_{2p+1},\ldots,d_{2p+2n})$. It is\vspace*{1.5pt} easy to verify
that the identity matrix $I_n$ generates a full $2^{2n} \times2n$
design. Therefore, the properties of $D$ depend on the matrix $V$ only.

For $s=\{c_1,c_2,\ldots,c_k\}$, a subset of $k \leq2n+2p$ columns of
$D$, define $j_k(s;D) = \sum_{i=1}^{2^{2n}} c_{s1} \cdots c_{sk}$, where
$c_{ij}$ is the $i${th} entry of $c_j$. The $j_k(s;D)$ values are
called the $J$-characteristics of design $D$ [\citet{DenTan99},
\citet{Tan01}]. It is evident that $|j_k(s;D)| \leq2^{2n}$. Following
Cheng, Li and Ye (\citeyear{CheLiYe04}), we define the aliasing index as $\rho_k(s) = \rho_k(s;D) =
 |j_k(s;D)|/2^{2n}$, which measures the amount of aliasing
among columns in $s$. It is obvious that $0 \leq\rho_k(s) \leq1$.
When $\rho_k(s) = 1$, the columns in $s$ are fully aliased with each
other and form a complete word of length $k$. It is equivalent to the
defining relations in regular designs. When $0 < \rho_k(s) < 1$, the
columns in $s$ are partially aliased with each other and form a partial
word of length $k$ with aliasing index $\rho_k(s)$. When $\rho_k(s) =
0$, the columns in $s$ are orthogonal and do not form a word.

Throughout this paper, for $\vec{i}$ to be a quaternary row vector, let
$f_{\vec{i}}$ be the number of times that $\vec{i}$ appears in the rows
of $V$. Define $\vec{w}=(w_1,\ldots,w_p)$ to be a word type that
describes the structure of a word. All $w_i$ are quaternary with the
following meanings. For $i=1,\ldots,p$, if $w_i=0$, none of the
$(2i-1)${th} and $(2i)${th} in $D$ are included in the word; if
$w_i=2$, both the $(2i-1)${th} and $(2i)${th} in $D$ are included in
the word; if $w_i$ is odd, either the $(2i-1)${th} or $(2i)${th} in
$D$ is included in the word. If there are $q$ odd entries in $\vec{w}$,
where $q<p$, there are $2^q$ different column choices. Therefore, we
denote $w_i=1$ or $3$ for different $i$ to represent different column
choices. For example, in $(1/16)${th}-fraction QC designs, there are
four possible forms of words, namely, $(1,1)$, $(1,3)$, $(3,1)$ and
$(3,3)$, representing the cases that select one column from the first
two columns of $D$ and select another column from the next two columns
of $D$.

Let $k_{\vec{w}}$ be the wordlength equation, or simply called
$k$-equation, of the word described by $\vec{w}$. In addition, denote
$C(p)$ by a $4^p \times p$ matrix consisting of all possible
combinations of quaternary entries. With reference to the matrix $V$,
$k_{\vec{w}}$ can be written as the linear combination of $f_{\vec
{i}}$, where $\vec{i}$ represents the $i${th} row of $C(p)$, that is,
$k_{\vec{w}} = \sum_{\vec{i} \in C(p)} c_{\vec{i}}f_{\vec{i}}$ for
$c_{\vec{i}}=0,1,2$. Furthermore, if there exists two $k$-equations
$k_{\vec{w}_1}$ and $k_{\vec{w}_2}$ with the corresponding coefficient
vectors $c_{\vec{i}}$ and $c'_{\vec{i}}$ in their summations, then we
define a code arithmetic (CA) operator $\oplus$ in the following way:
\[
k_{\vec{w}_1} \oplus k_{\vec{w}_2} = \biggl(\sum
_{\vec{i} \in C(p)} c_{\vec
{i}}f_{\vec{i}}\biggr) \oplus\biggl(
\sum_{\vec{i} \in C(p)} c'_{\vec{i}}f_{\vec
{i}}
\biggr) = \sum_{\vec{i} \in C(p)} L_w
\bigl(c_{\vec{i}}+c'_{\vec{i}}\bigr) f_{\vec
{i}} ,
\]
where $L_w(x)$ represents the Lee weight of $x$ and the Lee weights of
$0,1,2,3 \in Z_4$ are $0,1,2,1$, respectively. Notice that the
wordlength of a word is not equal to the value of $k$-equations directly,
but it is equal to that plus a constant showing the number of columns
among the first $2p$ columns of $D$ (generated from $V$) that are
included in the word.

The above definitions and concepts are demonstrated in the following example.

\begin{example}\label{ex1} 
Consider a general $(1/16)${th}-fraction QC design $D$ (i.e., $p=2$)
generated by a generator matrix $G=(V,I_n)$, where $V=(u,v)$ for
convenience. There are $16$ possible combinations of quaternary entries
for $\vec{i}=(i_1,i_2)$ for $i_1,i_2 \in\{0,1,2,3\}$. Given a word
formed by a specific group of columns $\vec{w}$, its $k$-equations
$k_{\vec{w}}$ can always be written as linear combinations of these 16
combinations of $\vec{i}$. For example,
\begin{eqnarray*}
k_{10}&= & 0(f_{00}+f_{01}+f_{02}+f_{03})\\
&&{}+1(f_{10}+f_{11}+f_{12}+f_{13}+f_{30}+f_{31} +f_{32}+f_{33})\\
&&{}+2(f_{20}+f_{21}+f_{22}+f_{23})
= l_1 , \\
k_{02}&= & 0(f_{00}+f_{02}+f_{10}+f_{12}+f_{20}+f_{22}+f_{30}+f_{32})\\
&&{}+2(f_{01}+f_{03} +f_{11}+f_{13}+f_{21}+f_{23}+f_{31}+f_{33})
= l_6,
\end{eqnarray*}
where $l_1$ and $l_6$ are defined in \citet{Zhaetal11}. If we perform
a CA operation on these two $k$-equations,
\begin{eqnarray*}
k_{10} \oplus k_{02} &= &
0(f_{00}+f_{02}+f_{21}+f_{23})\\
&&{}+1(f_{10}+f_{11}+f_{12}+f_{13} +f_{30}+f_{31}+f_{32}+f_{33})\\
&&{}+2(f_{01}+f_{03}+f_{20}+f_{22}).
\end{eqnarray*}
In the resulting $k$-equation, the coefficient of $f_{11}$ and $f_{21}$
come from $L_w(1+2)=1$ and $L_w(2+2)=0$, respectively.
\end{example}

For a simpler notation, we may write a set of $k$-equations into a matrix
form $K=CF$, where $K$ and $F$ are the $k$-equations and frequency
vectors, $C$ is the wordlength equation coefficient matrix or simply
called $k$-matrix. For $(1/4)${th}-fractions, $F=(f_0,f_1,f_2,f_3)^T$,
$K=(k_1,k_2)^T$ and the equations of $k_1$ and $k_2$ in \citet{PhoXu09} are rewritten as
\[
C=\pmatrix{ 0 & 1 & 2 & 1
\vspace*{2pt}\cr
0 & 2 & 0 & 2}.
\]
For $(1/16)${th}-fractions,
$F=(f_{00},f_{01},f_{02},f_{03},f_{10},f_{11},f_{12},f_{13},f_{20},f_{21},f_{22}, f_{23},\break f_{30},f_{31},f_{32},f_{33})^T$,
$K=(k_{01},k_{10},k_{02},k_{11},k_{13},k_{20},k_{12},k_{21},k_{22})^T$
and the equations of $l_1,\ldots,k_{10}$ in \citet{Zhaetal11} are
rewritten as
\[
C=\pmatrix{ 0 & 1 & 2 & 1 & 0 & 1 & 2
& 1 & 0 & 1 & 2 & 1 & 0 & 1 & 2 & 1
\vspace*{2pt}\cr
0 & 0 & 0 & 0 & 1 & 1 & 1 & 1 & 2 & 2 & 2 & 2 & 1 & 1 & 1 & 1
\vspace*{2pt}\cr
0 & 2 & 0 & 2 & 0 & 2 & 0 & 2 & 0 & 2 & 0 & 2 & 0 & 2 & 0 & 2
\vspace*{2pt}\cr
0 & 1 & 2 & 1 & 1 & 2 & 1 & 0 & 2 & 1 & 0 & 1 & 1 & 0 & 1 & 2
\vspace*{2pt}\cr
0 & 1 & 2 & 1 & 1 & 0 & 1 & 2 & 2 & 1 & 0 & 1 & 1 & 2 & 1 & 0
\vspace*{2pt}\cr
0 & 0 & 0 & 0 & 2 & 2 & 2 & 2 & 0 & 0 & 0 & 0 & 2 & 2 & 2 & 2
\vspace*{2pt}\cr
0 & 2 & 0 & 2 & 1 & 1 & 1 & 1 & 2 & 0 & 2 & 0 & 1 & 1 & 1 & 1
\vspace*{2pt}\cr
0 & 1 & 2 & 1 & 2 & 1 & 0 & 1 & 0 & 1 & 2 & 1 & 2 & 1 & 0 & 1
\vspace*{2pt}\cr
0 & 2 & 0 & 2 & 2 & 0 & 2 & 0 & 0 & 2 & 0 & 2 & 2 & 0 & 2 & 0}.
\]

The $k$-equations are ordered in the vector $K$ under the following
rules: (1) the position of $k_{\vec{i}_1}$ is on the front of that of
$k_{\vec{i}_2}$ if $\sum_{q=1}^p L_w(i_{1,q}) < \sum_{q=1}^p
L_w(i_{2,q})$; (2) if $\sum_{q=1}^p L_w(i_{1,q}) = \sum_{q=1}^p
L_w(i_{2,q})$, then the position of $k_{\vec{i}_1}$ is on the front of
that of $k_{\vec{i}_2}$ if $i_{1,u} < i_{2,u}$ and $i_{1,q} = i_{2,q}$
for all $0<q<u$, where $i_{1,q}$ and $i_{2,q}$ are the $q${th} entries
of $\vec{i}_1$ and $\vec{i}_2$, respectively. The frequency vector $F$
is ordered in the ascending order of its quaternary-coded decimal
counterpart. The $k$-matrix of higher-order-fraction QC designs ($p>1$)
will be discussed in the later part of this paper.

The aliasing index can be written in the form of $\rho=2^{-\lfloor
(a+\delta)/2 \rfloor}$, where $a$ is a linear combination of
frequencies. Therefore, we may write all $a$'s into a matrix form
$A=BF$, where $A$ is the aliasing index equation vector or simply
called a-equations, and $B$ is the aliasing index equation coefficient
matrix or simply called a-matrix. Generally speaking, the aliasing
index of each $k_{\vec{w}}$ is $\rho_{\vec{w} (\operatorname{mod} 2)}$, and $a_{\vec
{w} (\operatorname{mod} 2)}$ is a component of its order by definition. In addition,
$\delta=1$ if the sum of entries of $\vec{w}$ is even, or $0$
otherwise. According to \citet{PhoXu09}, for $(1/4)${th}-fractions,
there is only one aliasing index for $k_1$, so $A=(a_1)$ and $B=(0 1 0
1)$. For $(1/16)${th}-fractions in \citet{Zhaetal11}, there are
three aliasing indexes $A=(a_{01},a_{10},a_{11})$ and
\[
B=\pmatrix{ 0 & 1 & 0 & 1 & 0 & 1 & 0
& 1 & 0 & 1 & 0 & 1 & 0 & 1 & 0 & 1
\vspace*{2pt}\cr
0 & 0 & 0 & 0 & 1 & 1 & 1 & 1 & 0 & 0 & 0 & 0 & 1 & 1 & 1 & 1
\vspace*{2pt}\cr
0 & 1 & 0 & 1 & 1 & 0 & 1 & 0 & 0 & 1 & 0 & 1 & 1 & 0 & 1 & 0
\vspace*{2pt}\cr}.
\]
In general, the a-equations are ordered in the vector $A$ under similar
rules as $k$-equations in $K$.

\section{Some rules and theorems on the structure of quaternary-code designs}\label{sec3}

Given a general $k$-equation in $(1/4)^p$th-fraction QC designs $k_{\vec
{w}}=\sum_{\vec{i} \in C(p)} c_{\vec{i}}f_{\vec{i}}$, where $c_{\vec
{i}}=0,1$ or $2$, all entries of $\vec{w}$ are quaternary and $\vec{i}$
is the $i${th} row of $C(p)$,\vadjust{\goodbreak} we denote $\vec{w}=(\vec{w}_{l},\vec
{w}_{p-l})$ as a partition into two segments: the first segment has
length $l$ and the second segment has length $p-l$. Similarly, we
denote all $\vec{i}=(\vec{i}_{l},\vec{i}_{p-l})$. In addition, if $c$
is quaternary constant, $\vec{c}_{l}$ represents a vector of length $l$
that all entries are constant $c$.

The following rules suggest how a $k$-equation can be derived from
another $k$-equation. Rule~\ref{ru1} extends the $k$-equations in
$(1/4)^p$th-fraction QC designs to those in $(1/4)^{p+1}$th-fractions.

\begin{lemma}\label{ru1} 
Given a general $k$-equation in a $(1/4)^p$th-fraction QC design~$D_0$,\break
$k_{\vec{w}}=\sum_{\vec{i} \in C(p)} c_{\vec{i}}f_{\vec{i}}$. Then
all $k$-equations with $w_{l+1}=0$ in a $(1/4)^{p+1}$th-fraction QC
design $D$ can be expressed as $k_{(\vec{w}_{l},0,\vec{w}_{p-l})} = \sum_{s=0}^3
\sum_{\vec{i} \in C(p)} c_{\vec{i}}f_{(\vec{i}_{l},s,\vec{i}_{p-l})}$.
\end{lemma}

This result is obvious. If a $k$-equation consists of $w_{l+1}=0$, the
word described by this $k$-equation includes none of the $(2l+1)$th and
$(2l+2)$th columns of $D$. It acts like considering the same $k$-equation
in $(1/4)^p$th-fraction QC design $D_0$. This rule can be used to form
the basic $k$-equations for QC designs, which are stated in the following
corollaries.

\begin{corollary}\label{co1} 
For a general\vspace*{-1.5pt} $(1/4)^{p+1}$th-fraction QC design, $k_{(\vec{0}_l,1,\vec
{0}_{p-l})} = \sum_{\vec{i} \in C(p)} (f_{(\vec{i}_{l},1,\vec
{i}_{p-l})}+f_{(\vec{i}_{l},3,\vec{i}_{p-l})}+2 f_{(\vec{i}_{l},2,\vec
{i}_{p-l})})$, where $(\vec{i}_{l},\vec{i}_{p-l})$ represents the
$i${th} row of $C(p)$.
\end{corollary}

\begin{corollary}\label{co2} 
For a general\vspace*{-1.5pt} $(1/4)^{p+1}$th-fraction QC design, $k_{(\vec{0}_l,2,\vec
{0}_{p-l})} = 2\sum_{\vec{i} \in C(p)} (f_{(\vec{i}_{l},1,\vec
{i}_{p-l})}+f_{(\vec{i}_{l},3,\vec{i}_{p-l})})$, where $(\vec
{i}_{l},\vec{i}_{p-l})$ represents the $i${th} row of~$C(p)$.
\end{corollary}

The proofs of two corollaries are given in the last section. Rule~\ref{ru2}
considers the $k$-equations of a word that consists of only one out of
two binary columns generated from every quaternary column in $V$.

\begin{lemma}\label{ru2} 
Given a $k$-equation in a $(1/4)^p$th-fraction QC design $k_{\vec{1}_p}=\sum_{\vec{i} \in C(p)} c_{\vec{i}}f_{(\vec{i}_{p-1},i_p)}$, where
$i_p$ represents the last entry of $\vec{i}$, then $k_{(1,\vec{3}_p)}=\sum_{s=0}^3 \sum_{\vec{i} \in C(p)} c_{\vec{i}}f_{(s,\vec
{i}_{p-1},(i_p+s)\operatorname{mod}4)}$.
\end{lemma}

It provides a gateway to extend from the $k$-equations of
$(1/4)^p$th-fraction to $(1/4)^{p+1}$th-fraction, where the subscript
vectors of the $k$-equations are all odd entries. For examples, this rule
helps to extend from $k_1$ of $(1/4)${th}-fraction to $k_{13}$ of
$(1/16)${th}-fraction, or $k_{111}$ of $(1/64)${th}-fraction to
$k_{1333}$ of $(1/256)${th}-fraction.

Rule~\ref{ru3} provides a relationship between two $k$-equations of words with
slight difference in the columns chosen.

\begin{lemma}\label{ru3} 
Given a general $k$-equation in a $(1/4)^{p+1}$th-fraction QC design
$k_{\vec{w}}=k_{(\vec{w}_{l},s_1,\vec{w}_{p-l})}$, then $k_{(\vec
{w}_{l},s_2,\vec{w}_{p-l})}=k_{\vec{w}} \oplus k_{(\vec{0}_{l},2,\vec
{0}_{p-l})}$, where $s_1 = (s_2+2) \operatorname{mod}4$.\vadjust{\goodbreak}
\end{lemma}

The addition of $k_{(\vec{0}_l,2,\vec{0}_{p-l})}$ implies that a new
word is derived from the original word with additional inclusion of the
$(2l-1)${th} and $(2l)${th} columns from $D$, plus some columns in
$I_n$ so that the termwise multiplication of these additional columns
results in a vector of $1$, that is, a complete aliased structure.
Notice that the inclusion of a column twice is equivalent to the
exclusion of the column. In the case when $s_1$ is odd, the exchange
between $1$ and $3$ represents a derivation of different form of
$k$-equations when the word includes either the $(2l-1)${th} or
$(2l)${th} column only. On the other hand, when $s_1$ is even, the
exchange between $0$ and $2$ represents a derivation of the $k$-equation
of a new word that includes or excludes both the $(2l-1)${th} or
$(2l)${th} columns.

Let $C_z(p)$ be a subset of $C(p)$ for $z=0,1,2$ as follows. For $z$ is
even, $C_z(p) = \{ \vec{i} \in C(p)\dvtx i_1+\cdots+i_p = z (\operatorname{mod} 4) \}$;
otherwise, $C_1(p) = \{ \vec{i} \in C(p)\dvtx i_1+\cdots+i_p = 1$ or $3(\operatorname{mod} 4) \}$. Then the general structure of a $k$-equation, where all
entries of $\vec{w}$ are odd, can be derived in the following
theorem.\vspace*{-2pt}

\begin{theorem}\label{th1} 
In a $(1/4)^p$th-fraction QC design, for all odd entries of $\vec
{w}=\vec{1}_p$, a $k$-equation is expressed as $k_{\vec{w}} = 1 \sum_{\vec
{i} \in C_1(p)} f_{\vec{i}} + 2 \sum_{\vec{i} \in C_2(p)} f_{\vec{i}}$.\vspace*{-2pt}
\end{theorem}

There are $2^{2p-1}$ frequencies with coefficients $1$, $2^{2p-2}$
frequencies with coefficients $0$ and $2^{2p-2}$ frequencies with
coefficients $2$. Furthermore, among those $2^{2p-2}$ frequencies with
coefficients $2$, there are $2^{p-1}$ frequencies that all entries of
$\vec{i}_2$ are either $0$ or $2$. It is also the same for those
$2^{2p-2}$ frequencies with coefficients~$0$.\vspace*{-2pt}

\begin{example}\label{ex2} 
We consider a $k$-equation $k_{11}$ in a general $(1/16)${th}-fraction
QC design $D$. We can express $k_{11} =
1(f_{01}+f_{10}+f_{21}+f_{12}+f_{03}+f_{30}+f_{23}+f_{32})+2(f_{02}+f_{20}+f_{11}+f_{33})$,
that is, $C_0=\{(00),(22),(13),(31)\}$, $C_1=\{(01),(10),(21),(12),(03),(30),(23),(32)\}$ and $C_2=\{(02),(20),(11),(33)\}$. By counting
the above frequencies, there are $2^{2p-1}=8$ frequencies with
coefficient $1$, $2^{2p-2}=4$ frequencies with coefficients $0$ and
$2^{2p-2}=4$ frequencies with coefficients $2$. Furthermore, among
those four frequencies with coefficients~$2$, there are two frequencies
($f_{02}$ and $f_{20}$) that all entries of $\vec{i}_2$ are either $0$
or $2$. It is also the same for those frequencies with coefficients $0$
($f_{00}$ and $f_{22}$).\vspace*{-2pt}
\end{example}

The last rule defines the a-equation of a word accompanied with a
$k$-equation.\vspace*{-2pt}

\begin{lemma}\label{ru4} 
Given a general $k$-equation in a $(1/4)^p$th-fraction QC design $k_{\vec
{w}}$ as in Theorem~\ref{th1}, then the a-equation of the corresponding word is
$a_{\vec{w}}=a_{\vec{w} \operatorname{mod} 2}=\sum_{\vec{i} \in C_1(p)} f_{\vec{i}}$.\vspace*{-2pt}
\end{lemma}

Rule~\ref{ru4} implies that the aliasing index of a word depends only on the
number of odd entries in $\vec{w}$\vadjust{\goodbreak} and their positions, and the even
entries basically have no effects. For example, $k_{10}$ and $k_{12}$
are expected to share the same aliasing index $a_{10}$, but $k_{110}$
and $k_{011}$ are expected to have different aliasing indexes, the
prior has aliasing index $a_{110}$ and the latter has aliasing
index~$a_{011}$.

Among all $4^p$ $k$-equations for a general $(1/4)^p$th-fraction QC
design $D$, some of them are equivalent to others and some are
irrelevant. The following theorem considers these equivalences and
irrelevance and specifies a list of $k$-equations that are necessary to
be computed in order to obtain the properties of $D$.

\begin{theorem}\label{th2} 
Consider a general $(1/4)^p$th-fraction QC design $D$. There exists
$4^p$ possible combinations of $\vec{w}$ for $k$-equations. It is
necessary and sufficient to consider the following $\vec{w}$ in order
to obtain the properties of~$D$:
\begin{longlist}[(1)]
\item[(1)] $\vec{w}$ that all entries are even, except all entries are $0$;
\item[(2)] $\vec{w}$ that the first odd entry must be $1$ for $\vec{w}$
that consists of odd entries.
\end{longlist}
There are $2^p-1$ $k$-equations in the first group of $\vec{w}$ and
$2^{2p-1}-2^{p-1}$ $k$-equations in the second group.
\end{theorem}

\begin{example}\label{ex3} 
We consider a general $(1/16)${th}-fraction QC design $D$ and there
are 16 possible combinations of $\vec{w}$ listed in Example~\ref{ex1}.
According to Theorem~\ref{th2}, the first group of $\vec{w}$ has only even
entries. Except $\{0,0\}$, there are three combinations that satisfy
this situation and they are $\{0,2\}$, $\{2,0\}$ and $\{2,2\}$. For the
remaining 12 combinations (with at least one odd entry), these 6
combinations $\{0,3\}$, $\{2,3\}$, $\{3,0\}$, $\{3,1\}$, $\{3,2\}$, $\{
3,3\}$ are not included in consideration because the $k$-equations of
them are exactly equivalent to those with $\vec{w}=\{0,1\}$, $\{2,1\}
$, $\{1,0\}$, $\{1,3\}$, $\{1,2\}$, $\{1,1\}$, respectively. Therefore,
among all 16 possible combinations of $\vec{w}$, only 9 of them, 3 in
the first group and 6 in the second group, are necessary and sufficient
to be considered in order to determine the properties
of~$D$.
\end{example}

\section{Code arithmetic (CA) approach for generating wordlength
equations of $(1/64)${th}-fraction QC designs}\label{sec4}

This section extends the results of\break $(1/16)${th}-fraction QC designs
that appeared in \citet{Zhaetal11} and \citet{PhoMukXu12},
and sets of $k$-equations and a-equations for $(1/64)${th}-fractions QC
designs are generated using the theorems above. These \mbox{equations} are
applied to derive the design properties of $(1/64)${th}-fraction QC
designs.

Following Theorem~\ref{th2}, 35 $k$-equations are sufficient to determine the
properties of a $(1/64)${th}-fraction QC design. Specifically, seven
of them belong to the first group and 28 of them belong to the second
group. Using the CA approach, we derive these 35 $k$-equations and their
corresponding a-equations from the $k$-equations of $(1/4)${th}- and
$(1/16)${th}-fractions QC designs. First, we define $C(2) $ to be a
$16 \times2$ matrix consisting of all 16 possible combinations of
quaternary entries. Throughout this section, we express all $k$-equations
as a row in the $k$-matrix for clear and convenient notation.

Rule~\ref{ru1} and two corollaries are applied to obtain $k$-equations where $\vec
{i}$ contains at least one $0$. More explicitly, to obtain $k$-equations
with two $0$s in $\vec{i}$, that is, $k_{100}$, $k_{010}$, $k_{001}$,
$k_{200}$, $k_{020}$ and $k_{002}$, we apply Corollaries~\ref{co1} and~\ref{co2} with
$l=0,1,2$. For example, for $k_{100}$, we apply Corollary~\ref{co1} with $l=0$.
This yields a $k$-equation where, for $j,k=0,1,2,3$, the coefficients of
$f_{0jk}$, $f_{1jk}$, $f_{2jk}$ and $f_{3jk}$ are $0$, $1$, $2$ and
$1$, respectively. Rule~\ref{ru4} suggests $a_{100}$, the a-equations of
$k_{100}$, such that the coefficients of $f_{0jk}$, $f_{1jk}$,
$f_{2jk}$ and $f_{3jk}$ are $0$, $1$, $0$ and $1$,
respectively.

For all $k$-equations with one $0$ in $\vec{i}$, we consider applying
Rule~\ref{ru1} on $k_{11}$, $k_{13}$, $k_{12}$, $k_{21}$ and $k_{22}$ with
different $l$. This leads to $k_{011}$, $k_{013}$, $k_{012}$,
$k_{021}$, $k_{022}$ when $l=0$, $k_{101}$, $k_{103}$, $k_{102}$,
$k_{201}$, $k_{202}$ when $l=1$ and $k_{110}$, $k_{130}$, $k_{120}$,
$k_{210}$, $k_{220}$ when $l=2$. For example, for $k_{101}$, Rule~\ref{ru1}
suggests that $\vec{w}_l=\vec{w}_{p-l}=1$. Then for every row of
$C(2)$, denoted as $(c_1,c_2)$, the coefficients of
$f_{(c_1,0,c_2)}$,$f_{(c_1,1,c_2)}$,$f_{(c_1,2,c_2)}$,$f_{(c_1,3,c_2)}$
in $k_{101}$ are all equal to the coefficient of $f_{(c_1,c_2)}$ in $k_{11}$.

It is straightforward to substitute $0$s in all $k$-equations mentioned
above with $2$ by Rule~\ref{ru3}. By changing one $0$ into $2$ in $\vec{i}$, we
obtain $k_{102}$, $k_{120}$, $k_{012}$, $k_{210}$, $k_{021}$,
$k_{201}$, $k_{202}$, $k_{220}$, $k_{022}$, $k_{112}$, $k_{132}$,
$k_{122}$, $k_{212}$, $k_{222}$, $k_{121}$, $k_{123}$, $k_{221}$,
$k_{211}$ and $k_{213}$. For example, in order to obtain $k_{121}$,
Rule~\ref{ru3} suggests performing a CA operation $k_{121} = k_{101} \oplus
k_{020}$. The a-equation of $k_{121}$ is equal to
$a_{101}$.

Rule~\ref{ru2} is applied in order to obtain the $k$-equations with all odd
entries in $\vec{i}$, including $k_{111}$, $k_{113}$, $k_{131}$ and~$k_{133}$.
According to Rule~\ref{ru2}, $k_{133}$ can be derived from~$k_{11}$.
For every row of $C(2)$, the first and second entries are considered as
$\vec{i}_{p-1}$ and $i_p$, respectively. For example, $\vec{i}_{p-1}=1$
and $i_p=0$ for $f_{10}$. Then we can determine the coefficients of
frequency vectors in $k_{133}$ from those in $k_{11}$. Consider $\vec
{i}=(10)$, for example. The coefficient of $f_{10}$ in $k_{11}$ is $1$.
This implies $f_{010}=f_{111}=f_{212}=f_{313}=1$ in
$k_{133}$ for $s=0,1,2,3$. Consider $\vec{i}=(02)$ as another example.
The coefficient of $f_{02}$ in $k_{11}$ is $2$. This implies $f_{002}=f_{103}=f_{200}=f_{301}=2$ in $k_{133}$ for $s=0,1,2,3$. The
other three $k$-equations without $0$s in $\vec{i}$ can be derived from
$k_{133}$ via the CA operations suggested in Rule~\ref{ru3}: $k_{111} =
(k_{133} \oplus k_{020}) \oplus k_{002}$, $k_{113} = k_{133} \oplus
k_{020}$, and $k_{131} = k_{133} \oplus k_{002}$. The a-equations of
$k_{111}$, $k_{113}$, $k_{131}$ and $k_{133}$ are the
same.\looseness=-1

There are in total 35 $k$-equations and 7 a-equations in $K$ and $A$,
respectively. Similar to $(1/16)${th}-fraction QC designs, we may
rewrite these $k$-equations and a-equations into matrix forms where
\begin{eqnarray*}
K&= & (k_{001},k_{010},k_{100},k_{002},k_{011},k_{013},k_{020},k_{101},k_{103},k_{110},k_{130},k_{200},k_{012},
\\
&&\hspace*{4pt}{} k_{021},k_{102},k_{111},k_{113},k_{131},k_{133},k_{120},k_{201},k_{210},k_{022},k_{112},k_{132},k_{121},
\\
&&\hspace*{100pt}{}
k_{123},k_{202},k_{211},k_{213},k_{220},k_{122},k_{212},k_{221},k_{222})^T,
\end{eqnarray*}
\begin{eqnarray*}
F&= & (f_{000},f_{001},f_{002},f_{003},f_{010},f_{011},f_{012},f_{013},f_{020},f_{021},f_{022},f_{023},f_{030},
\\
&&\hspace*{4pt}{} f_{031},f_{032},f_{033},f_{100},f_{101},f_{102},f_{103},f_{110},f_{111},f_{112},f_{113},f_{120},f_{121},
\\
&&\hspace*{4pt}{}f_{122},f_{123},f_{130},f_{131},f_{132},f_{133},f_{200},f_{201},f_{202},f_{203},f_{210},f_{211},f_{212},
\\
&&\hspace*{4pt}{} f_{213},f_{220},f_{221},f_{222},f_{223},f_{230},f_{231},f_{232},f_{233},f_{300},f_{301},f_{302},f_{303},
\\
&&\hspace*{32pt}{}
f_{310},f_{311},f_{312},f_{313},f_{320},f_{321},f_{322},f_{323},f_{330},f_{331},f_{332},f_{333})^T,
\end{eqnarray*}
\[
{\fontsize{10}{12}{\selectfont
C=
\lleft( %
\begin{array} {@{}c@{}c@{}c@{}c@{}c@{}c@{}c@{}c@{}c@{}c@{}c@{}c@{}c@{}c@{}c@{}c@{}c@{}c@{}c@{}c@{}c@{}c@{}c@{}c@{}
c@{}c@{}c@{}c@{}c@{}c@{}c@{}c@{}c@{}c@{}c@{}c@{}c@{}c@{}c@{}c@{}c@{}c@{}c@{}c@{}c@{}c@{}c@{}c@{}c@{}c@{}c@{}c@{}
c@{}c@{}c@{}c@{}c@{}c@{}c@{}c@{}c@{}c@{}c@{}c@{}}
0 & 1 & 2 & 1 & 0 & 1 & 2 & 1 &
0 & 1 & 2 & 1 & 0 & 1 & 2 & 1 & 0 & 1 & 2 & 1 & 0 & 1 & 2 & 1 & 0 & 1 & 2 & 1 & 0
& 1 & 2 & 1 & 0 & 1 & 2 & 1 & 0 & 1 & 2 & 1 & 0 & 1 & 2 & 1 & 0 & 1 & 2 & 1 & 0 &
1 & 2 & 1 & 0 & 1 & 2 & 1 & 0 & 1 & 2 & 1 & 0 & 1 & 2 & 1
\\
0 & 0 & 0 & 0 & 1 & 1 & 1 & 1 & 2 & 2 & 2 & 2 & 1 & 1 & 1 & 1 & 0 & 0 & 0 & 0 &
1 & 1 & 1 & 1 & 2 & 2 & 2 & 2 & 1 & 1 & 1 & 1 & 0 & 0 & 0 & 0 & 1 & 1 & 1 & 1 & 2
& 2 & 2 & 2 & 1 & 1 & 1 & 1 & 0 & 0 & 0 & 0 & 1 & 1 & 1 & 1 & 2 & 2 & 2 & 2 & 1 &
1 & 1 & 1
\\
0 & 0 & 0 & 0 & 0 & 0 & 0 & 0 & 0 & 0 & 0 & 0 & 0 & 0 & 0 & 0 & 1 & 1 & 1 & 1 &
1 & 1 & 1 & 1 & 1 & 1 & 1 & 1 & 1 & 1 & 1 & 1 & 2 & 2 & 2 & 2 & 2 & 2 & 2 & 2 & 2
& 2 & 2 & 2 & 2 & 2 & 2 & 2 & 1 & 1 & 1 & 1 & 1 & 1 & 1 & 1 & 1 & 1 & 1 & 1 & 1 &
1 & 1 & 1
\\
0 & 2 & 0 & 2 & 0 & 2 & 0 & 2 & 0 & 2 & 0 & 2 & 0 & 2 & 0 & 2 & 0 & 2 & 0 & 2 &
0 & 2 & 0 & 2 & 0 & 2 & 0 & 2 & 0 & 2 & 0 & 2 & 0 & 2 & 0 & 2 & 0 & 2 & 0 & 2 & 0
& 2 & 0 & 2 & 0 & 2 & 0 & 2 & 0 & 2 & 0 & 2 & 0 & 2 & 0 & 2 & 0 & 2 & 0 & 2 & 0 &
2 & 0 & 2
\\
0 & 1 & 2 & 1 & 1 & 2 & 1 & 0 & 2 & 1 & 0 & 1 & 1 & 0 & 1 & 2 & 0 & 1 & 2 & 1 &
1 & 2 & 1 & 0 & 2 & 1 & 0 & 1 & 1 & 0 & 1 & 2 & 0 & 1 & 2 & 1 & 1 & 2 & 1 & 0 & 2
& 1 & 0 & 1 & 1 & 0 & 1 & 2 & 0 & 1 & 2 & 1 & 1 & 2 & 1 & 0 & 2 & 1 & 0 & 1 & 1 &
0 & 1 & 2
\\
0 & 1 & 2 & 1 & 1 & 0 & 1 & 2 & 2 & 1 & 0 & 1 & 1 & 2 & 1 & 0 & 0 & 1 & 2 & 1 &
1 & 0 & 1 & 2 & 2 & 1 & 0 & 1 & 1 & 2 & 1 & 0 & 0 & 1 & 2 & 1 & 1 & 0 & 1 & 2 & 2
& 1 & 0 & 1 & 1 & 2 & 1 & 0 & 0 & 1 & 2 & 1 & 1 & 0 & 1 & 2 & 2 & 1 & 0 & 1 & 1 &
2 & 1 & 0
\\
0 & 0 & 0 & 0 & 2 & 2 & 2 & 2 & 0 & 0 & 0 & 0 & 2 & 2 & 2 & 2 & 0 & 0 & 0 & 0 &
2 & 2 & 2 & 2 & 0 & 0 & 0 & 0 & 2 & 2 & 2 & 2 & 0 & 0 & 0 & 0 & 2 & 2 & 2 & 2 & 0
& 0 & 0 & 0 & 2 & 2 & 2 & 2 & 0 & 0 & 0 & 0 & 2 & 2 & 2 & 2 & 0 & 0 & 0 & 0 & 2 &
2 & 2 & 2
\\
0 & 1 & 2 & 1 & 0 & 1 & 2 & 1 & 0 & 1 & 2 & 1 & 0 & 1 & 2 & 1 & 1 & 2 & 1 & 0 &
1 & 2 & 1 & 0 & 1 & 2 & 1 & 0 & 1 & 2 & 1 & 0 & 2 & 1 & 0 & 1 & 2 & 1 & 0 & 1 & 2
& 1 & 0 & 1 & 2 & 1 & 0 & 1 & 1 & 0 & 1 & 2 & 1 & 0 & 1 & 2 & 1 & 0 & 1 & 2 & 1 &
0 & 1 & 2
\\
0 & 1 & 2 & 1 & 0 & 1 & 2 & 1 & 0 & 1 & 2 & 1 & 0 & 1 & 2 & 1 & 1 & 0 & 1 & 2 &
1 & 0 & 1 & 2 & 1 & 0 & 1 & 2 & 1 & 0 & 1 & 2 & 2 & 1 & 0 & 1 & 2 & 1 & 0 & 1 & 2
& 1 & 0 & 1 & 2 & 1 & 0 & 1 & 1 & 2 & 1 & 0 & 1 & 2 & 1 & 0 & 1 & 2 & 1 & 0 & 1 &
2 & 1 & 0
\\
0 & 0 & 0 & 0 & 1 & 1 & 1 & 1 & 2 & 2 & 2 & 2 & 1 & 1 & 1 & 1 & 1 & 1 & 1 & 1 &
2 & 2 & 2 & 2 & 1 & 1 & 1 & 1 & 0 & 0 & 0 & 0 & 2 & 2 & 2 & 2 & 1 & 1 & 1 & 1 & 0
& 0 & 0 & 0 & 1 & 1 & 1 & 1 & 1 & 1 & 1 & 1 & 0 & 0 & 0 & 0 & 1 & 1 & 1 & 1 & 2 &
2 & 2 & 2
\\
0 & 0 & 0 & 0 & 1 & 1 & 1 & 1 & 2 & 2 & 2 & 2 & 1 & 1 & 1 & 1 & 1 & 1 & 1 & 1 &
0 & 0 & 0 & 0 & 1 & 1 & 1 & 1 & 2 & 2 & 2 & 2 & 2 & 2 & 2 & 2 & 1 & 1 & 1 & 1 & 0
& 0 & 0 & 0 & 1 & 1 & 1 & 1 & 1 & 1 & 1 & 1 & 2 & 2 & 2 & 2 & 1 & 1 & 1 & 1 & 0 &
0 & 0 & 0
\\
0 & 0 & 0 & 0 & 0 & 0 & 0 & 0 & 0 & 0 & 0 & 0 & 0 & 0 & 0 & 0 & 2 & 2 & 2 & 2 &
2 & 2 & 2 & 2 & 2 & 2 & 2 & 2 & 2 & 2 & 2 & 2 & 0 & 0 & 0 & 0 & 0 & 0 & 0 & 0 & 0
& 0 & 0 & 0 & 0 & 0 & 0 & 0 & 2 & 2 & 2 & 2 & 2 & 2 & 2 & 2 & 2 & 2 & 2 & 2 & 2 &
2 & 2 & 2
\\
0 & 2 & 0 & 2 & 1 & 1 & 1 & 1 & 2 & 0 & 2 & 0 & 1 & 1 & 1 & 1 & 0 & 2 & 0 & 2 &
1 & 1 & 1 & 1 & 2 & 0 & 2 & 0 & 1 & 1 & 1 & 1 & 0 & 2 & 0 & 2 & 1 & 1 & 1 & 1 & 2
& 0 & 2 & 0 & 1 & 1 & 1 & 1 & 0 & 2 & 0 & 2 & 1 & 1 & 1 & 1 & 2 & 0 & 2 & 0 & 1 &
1 & 1 & 1
\\
0 & 1 & 2 & 1 & 2 & 1 & 0 & 1 & 0 & 1 & 2 & 1 & 2 & 1 & 0 & 1 & 0 & 1 & 2 & 1 &
2 & 1 & 0 & 1 & 0 & 1 & 2 & 1 & 2 & 1 & 0 & 1 & 0 & 1 & 2 & 1 & 2 & 1 & 0 & 1 & 0
& 1 & 2 & 1 & 2 & 1 & 0 & 1 & 0 & 1 & 2 & 1 & 2 & 1 & 0 & 1 & 0 & 1 & 2 & 1 & 2 &
1 & 0 & 1
\\
0 & 2 & 0 & 2 & 0 & 2 & 0 & 2 & 0 & 2 & 0 & 2 & 0 & 2 & 0 & 2 & 1 & 1 & 1 & 1 &
1 & 1 & 1 & 1 & 1 & 1 & 1 & 1 & 1 & 1 & 1 & 1 & 2 & 0 & 2 & 0 & 2 & 0 & 2 & 0 & 2
& 0 & 2 & 0 & 2 & 0 & 2 & 0 & 1 & 1 & 1 & 1 & 1 & 1 & 1 & 1 & 1 & 1 & 1 & 1 & 1 &
1 & 1 & 1
\\
0 & 1 & 2 & 1 & 1 & 2 & 1 & 0 & 2 & 1 & 0 & 1 & 1 & 0 & 1 & 2 & 1 & 2 & 1 & 0 &
2 & 1 & 0 & 1 & 1 & 0 & 1 & 2 & 0 & 1 & 2 & 1 & 2 & 1 & 0 & 1 & 1 & 0 & 1 & 2 & 0
& 1 & 2 & 1 & 1 & 2 & 1 & 0 & 1 & 0 & 1 & 2 & 0 & 1 & 2 & 1 & 1 & 2 & 1 & 0 & 2 &
1 & 0 & 1
\\
0 & 1 & 2 & 1 & 1 & 0 & 1 & 2 & 2 & 1 & 0 & 1 & 1 & 2 & 1 & 0 & 1 & 0 & 1 & 2 &
2 & 1 & 0 & 1 & 1 & 2 & 1 & 0 & 0 & 1 & 2 & 1 & 2 & 1 & 0 & 1 & 1 & 2 & 1 & 0 & 0
& 1 & 2 & 1 & 1 & 0 & 1 & 2 & 1 & 2 & 1 & 0 & 0 & 1 & 2 & 1 & 1 & 0 & 1 & 2 & 2 &
1 & 0 & 1
\\
0 & 1 & 2 & 1 & 1 & 0 & 1 & 2 & 2 & 1 & 0 & 1 & 1 & 2 & 1 & 0 & 1 & 2 & 1 & 0 &
0 & 1 & 2 & 1 & 1 & 0 & 1 & 2 & 2 & 1 & 0 & 1 & 2 & 1 & 0 & 1 & 1 & 2 & 1 & 0 & 0
& 1 & 2 & 1 & 1 & 0 & 1 & 2 & 1 & 0 & 1 & 2 & 2 & 1 & 0 & 1 & 1 & 2 & 1 & 0 & 0 &
1 & 2 & 1
\\
0 & 1 & 2 & 1 & 1 & 2 & 1 & 0 & 2 & 1 & 0 & 1 & 1 & 0 & 1 & 2 & 1 & 0 & 1 & 2 &
0 & 1 & 2 & 1 & 1 & 2 & 1 & 0 & 2 & 1 & 0 & 1 & 2 & 1 & 0 & 1 & 1 & 0 & 1 & 2 & 0
& 1 & 2 & 1 & 1 & 2 & 1 & 0 & 1 & 2 & 1 & 0 & 2 & 1 & 0 & 1 & 1 & 0 & 1 & 2 & 0 &
1 & 2 & 1
\\
0 & 0 & 0 & 0 & 2 & 2 & 2 & 2 & 0 & 0 & 0 & 0 & 2 & 2 & 2 & 2 & 1 & 1 & 1 & 1 &
1 & 1 & 1 & 1 & 1 & 1 & 1 & 1 & 1 & 1 & 1 & 1 & 2 & 2 & 2 & 2 & 0 & 0 & 0 & 0 & 2
& 2 & 2 & 2 & 0 & 0 & 0 & 0 & 1 & 1 & 1 & 1 & 1 & 1 & 1 & 1 & 1 & 1 & 1 & 1 & 1 &
1 & 1 & 1
\\
0 & 1 & 2 & 1 & 0 & 1 & 2 & 1 & 0 & 1 & 2 & 1 & 0 & 1 & 2 & 1 & 2 & 1 & 0 & 1 &
2 & 1 & 0 & 1 & 2 & 1 & 0 & 1 & 2 & 1 & 0 & 1 & 0 & 1 & 2 & 1 & 0 & 1 & 2 & 1 & 0
& 1 & 2 & 1 & 0 & 1 & 2 & 1 & 2 & 1 & 0 & 1 & 2 & 1 & 0 & 1 & 2 & 1 & 0 & 1 & 2 &
1 & 0 & 1
\\
0 & 0 & 0 & 0 & 1 & 1 & 1 & 1 & 2 & 2 & 2 & 2 & 1 & 1 & 1 & 1 & 2 & 2 & 2 & 2 &
1 & 1 & 1 & 1 & 0 & 0 & 0 & 0 & 1 & 1 & 1 & 1 & 0 & 0 & 0 & 0 & 1 & 1 & 1 & 1 & 2
& 2 & 2 & 2 & 1 & 1 & 1 & 1 & 2 & 2 & 2 & 2 & 1 & 1 & 1 & 1 & 0 & 0 & 0 & 0 & 1 &
1 & 1 & 1
\\
0 & 2 & 0 & 2 & 2 & 0 & 2 & 0 & 0 & 2 & 0 & 2 & 2 & 0 & 2 & 0 & 0 & 2 & 0 & 2 &
2 & 0 & 2 & 0 & 0 & 2 & 0 & 2 & 2 & 0 & 2 & 0 & 0 & 2 & 0 & 2 & 2 & 0 & 2 & 0 & 0
& 2 & 0 & 2 & 2 & 0 & 2 & 0 & 0 & 2 & 0 & 2 & 2 & 0 & 2 & 0 & 0 & 2 & 0 & 2 & 2 &
0 & 2 & 0
\\
0 & 2 & 0 & 2 & 1 & 1 & 1 & 1 & 2 & 0 & 2 & 0 & 1 & 1 & 1 & 1 & 1 & 1 & 1 & 1 &
2 & 0 & 2 & 0 & 1 & 1 & 1 & 1 & 0 & 2 & 0 & 2 & 2 & 0 & 2 & 0 & 1 & 1 & 1 & 1 & 0
& 2 & 0 & 2 & 1 & 1 & 1 & 1 & 1 & 1 & 1 & 1 & 0 & 2 & 0 & 2 & 1 & 1 & 1 & 1 & 2 &
0 & 2 & 0
\\
0 & 2 & 0 & 2 & 1 & 1 & 1 & 1 & 2 & 0 & 2 & 0 & 1 & 1 & 1 & 1 & 1 & 1 & 1 & 1 &
0 & 2 & 0 & 2 & 1 & 1 & 1 & 1 & 2 & 0 & 2 & 0 & 2 & 0 & 2 & 0 & 1 & 1 & 1 & 1 & 0
& 2 & 0 & 2 & 1 & 1 & 1 & 1 & 1 & 1 & 1 & 1 & 2 & 0 & 2 & 0 & 1 & 1 & 1 & 1 & 0 &
2 & 0 & 2
\\
0 & 1 & 2 & 1 & 2 & 1 & 0 & 1 & 0 & 1 & 2 & 1 & 2 & 1 & 0 & 1 & 1 & 2 & 1 & 0 &
1 & 0 & 1 & 2 & 1 & 2 & 1 & 0 & 1 & 0 & 1 & 2 & 2 & 1 & 0 & 1 & 0 & 1 & 2 & 1 & 2
& 1 & 0 & 1 & 0 & 1 & 2 & 1 & 1 & 0 & 1 & 2 & 1 & 2 & 1 & 0 & 1 & 0 & 1 & 2 & 1 &
2 & 1 & 0
\\
0 & 1 & 2 & 1 & 2 & 1 & 0 & 1 & 0 & 1 & 2 & 1 & 2 & 1 & 0 & 1 & 1 & 0 & 1 & 2 &
1 & 2 & 1 & 0 & 1 & 0 & 1 & 2 & 1 & 2 & 1 & 0 & 2 & 1 & 0 & 1 & 0 & 1 & 2 & 1 & 2
& 1 & 0 & 1 & 0 & 1 & 2 & 1 & 1 & 2 & 1 & 0 & 1 & 0 & 1 & 2 & 1 & 2 & 1 & 0 & 1 &
0 & 1 & 2
\\
0 & 2 & 0 & 2 & 0 & 2 & 0 & 2 & 0 & 2 & 0 & 2 & 0 & 2 & 0 & 2 & 2 & 0 & 2 & 0 &
2 & 0 & 2 & 0 & 2 & 0 & 2 & 0 & 2 & 0 & 2 & 0 & 0 & 2 & 0 & 2 & 0 & 2 & 0 & 2 & 0
& 2 & 0 & 2 & 0 & 2 & 0 & 2 & 2 & 0 & 2 & 0 & 2 & 0 & 2 & 0 & 2 & 0 & 2 & 0 & 2 &
0 & 2 & 0
\\
0 & 1 & 2 & 1 & 1 & 2 & 1 & 0 & 2 & 1 & 0 & 1 & 1 & 0 & 1 & 2 & 2 & 1 & 0 & 1 &
1 & 0 & 1 & 2 & 0 & 1 & 2 & 1 & 1 & 2 & 1 & 0 & 0 & 1 & 2 & 1 & 1 & 2 & 1 & 0 & 2
& 1 & 0 & 1 & 1 & 0 & 1 & 2 & 2 & 1 & 0 & 1 & 1 & 0 & 1 & 2 & 0 & 1 & 2 & 1 & 1 &
2 & 1 & 0
\\
0 & 1 & 2 & 1 & 1 & 0 & 1 & 2 & 2 & 1 & 0 & 1 & 1 & 2 & 1 & 0 & 2 & 1 & 0 & 1 &
1 & 2 & 1 & 0 & 0 & 1 & 2 & 1 & 1 & 0 & 1 & 2 & 0 & 1 & 2 & 1 & 1 & 0 & 1 & 2 & 2
& 1 & 0 & 1 & 1 & 2 & 1 & 0 & 2 & 1 & 0 & 1 & 1 & 2 & 1 & 0 & 0 & 1 & 2 & 1 & 1 &
0 & 1 & 2
\\
0 & 0 & 0 & 0 & 2 & 2 & 2 & 2 & 0 & 0 & 0 & 0 & 2 & 2 & 2 & 2 & 2 & 2 & 2 & 2 &
0 & 0 & 0 & 0 & 2 & 2 & 2 & 2 & 0 & 0 & 0 & 0 & 0 & 0 & 0 & 0 & 2 & 2 & 2 & 2 & 0
& 0 & 0 & 0 & 2 & 2 & 2 & 2 & 2 & 2 & 2 & 2 & 0 & 0 & 0 & 0 & 2 & 2 & 2 & 2 & 0 &
0 & 0 & 0
\\
0 & 2 & 0 & 2 & 2 & 0 & 2 & 0 & 0 & 2 & 0 & 2 & 2 & 0 & 2 & 0 & 1 & 1 & 1 & 1 &
1 & 1 & 1 & 1 & 1 & 1 & 1 & 1 & 1 & 1 & 1 & 1 & 2 & 0 & 2 & 0 & 0 & 2 & 0 & 2 & 2
& 0 & 2 & 0 & 0 & 2 & 0 & 2 & 1 & 1 & 1 & 1 & 1 & 1 & 1 & 1 & 1 & 1 & 1 & 1 & 1 &
1 & 1 & 1
\\
0 & 2 & 0 & 2 & 1 & 1 & 1 & 1 & 2 & 0 & 2 & 0 & 1 & 1 & 1 & 1 & 2 & 0 & 2 & 0 &
1 & 1 & 1 & 1 & 0 & 2 & 0 & 2 & 1 & 1 & 1 & 1 & 0 & 2 & 0 & 2 & 1 & 1 & 1 & 1 & 2
& 0 & 2 & 0 & 1 & 1 & 1 & 1 & 2 & 0 & 2 & 0 & 1 & 1 & 1 & 1 & 0 & 2 & 0 & 2 & 1 &
1 & 1 & 1
\\
0 & 1 & 2 & 1 & 2 & 1 & 0 & 1 & 0 & 1 & 2 & 1 & 2 & 1 & 0 & 1 & 2 & 1 & 0 & 1 &
0 & 1 & 2 & 1 & 2 & 1 & 0 & 1 & 0 & 1 & 2 & 1 & 0 & 1 & 2 & 1 & 2 & 1 & 0 & 1 & 0
& 1 & 2 & 1 & 2 & 1 & 0 & 1 & 2 & 1 & 0 & 1 & 0 & 1 & 2 & 1 & 2 & 1 & 0 & 1 & 0 &
1 & 2 & 1
\\
0 & 2 & 0 & 2 & 2 & 0 & 2 & 0 & 0 & 2 & 0 & 2 & 2 & 0 & 2 & 0 & 2 & 0 & 2 & 0 &
0 & 2 & 0 & 2 & 2 & 0 & 2 & 0 & 0 & 2 & 0 & 2 & 0 & 2 & 0 & 2 & 2 & 0 & 2 & 0 & 0
& 2 & 0 & 2 & 2 & 0 & 2 & 0 & 2 & 0 & 2 & 0 & 0 & 2 & 0 & 2 & 2 & 0 & 2 & 0 & 0 &
2 & 0 & 2
\end{array} %
\rright)}},
\]
\[
A=(a_{001},a_{010},a_{100},a_{011},a_{101},a_{110},a_{111}),
\]
\[
{\fontsize{9.9}{12}{\selectfont
B=\lleft( %
\begin{array} {@{}c@{}c@{}c@{}c@{}c@{}c@{}c@{}c@{}c@{}c@{}c@{}c@{}c@{}c@{}c@{}c@{}c@{}c@{}c@{}c@{}c@{}c@{}c@{}c@{}c@{}c@{}c@{}c@{}c@{}c@{}c@{}c@{}c@{}c@{}c@{}c@{}c@{}c@{}c@{}c@{}c@{}c@{}c@{}c@{}c@{}c@{}c@{}c@{}c@{}c@{}c@{}c@{}c@{}c@{}c@{}c@{}c@{}c@{}c@{}c@{}c@{}c@{}c@{}c@{}}
0 & 1 & 0 & 1 & 0 & 1 & 0 & 1 & 0 & 1 & 0 & 1 & 0 & 1 & 0 & 1 & 0 & 1 & 0 & 1 & 0
& 1 & 0 & 1 & 0 & 1 & 0 & 1 & 0 & 1 & 0 & 1 & 0 & 1 & 0 & 1 & 0 & 1 & 0 & 1 & 0 &
1 & 0 & 1 & 0 & 1 & 0 & 1 & 0 & 1 & 0 & 1 & 0 & 1 & 0 & 1 & 0 & 1 & 0 & 1 & 0 & 1
& 0 & 1
\\
0 & 0 & 0 & 0 & 1 & 1 & 1 & 1 & 0 & 0 & 0 & 0 & 1 & 1 & 1 & 1 & 0 & 0 & 0 & 0 &
1 & 1 & 1 & 1 & 0 & 0 & 0 & 0 & 1 & 1 & 1 & 1 & 0 & 0 & 0 & 0 & 1 & 1 & 1 & 1 & 0
& 0 & 0 & 0 & 1 & 1 & 1 & 1 & 0 & 0 & 0 & 0 & 1 & 1 & 1 & 1 & 0 & 0 & 0 & 0 & 1 &
1 & 1 & 1
\\
0 & 0 & 0 & 0 & 0 & 0 & 0 & 0 & 0 & 0 & 0 & 0 & 0 & 0 & 0 & 0 & 1 & 1 & 1 & 1 &
1 & 1 & 1 & 1 & 1 & 1 & 1 & 1 & 1 & 1 & 1 & 1 & 0 & 0 & 0 & 0 & 0 & 0 & 0 & 0 & 0
& 0 & 0 & 0 & 0 & 0 & 0 & 0 & 1 & 1 & 1 & 1 & 1 & 1 & 1 & 1 & 1 & 1 & 1 & 1 & 1 &
1 & 1 & 1
\\
0 & 1 & 0 & 1 & 1 & 0 & 1 & 0 & 0 & 1 & 0 & 1 & 1 & 0 & 1 & 0 & 0 & 1 & 0 & 1 &
1 & 0 & 1 & 0 & 0 & 1 & 0 & 1 & 1 & 0 & 1 & 0 & 0 & 1 & 0 & 1 & 1 & 0 & 1 & 0 & 0
& 1 & 0 & 1 & 1 & 0 & 1 & 0 & 0 & 1 & 0 & 1 & 1 & 0 & 1 & 0 & 0 & 1 & 0 & 1 & 1 &
0 & 1 & 0
\\
0 & 1 & 0 & 1 & 0 & 1 & 0 & 1 & 0 & 1 & 0 & 1 & 0 & 1 & 0 & 1 & 1 & 0 & 1 & 0 &
1 & 0 & 1 & 0 & 1 & 0 & 1 & 0 & 1 & 0 & 1 & 0 & 0 & 1 & 0 & 1 & 0 & 1 & 0 & 1 & 0
& 1 & 0 & 1 & 0 & 1 & 0 & 1 & 1 & 0 & 1 & 0 & 1 & 0 & 1 & 0 & 1 & 0 & 1 & 0 & 1 &
0 & 1 & 0
\\
0 & 0 & 0 & 0 & 1 & 1 & 1 & 1 & 0 & 0 & 0 & 0 & 1 & 1 & 1 & 1 & 1 & 1 & 1 & 1 &
0 & 0 & 0 & 0 & 1 & 1 & 1 & 1 & 0 & 0 & 0 & 0 & 0 & 0 & 0 & 0 & 1 & 1 & 1 & 1 & 0
& 0 & 0 & 0 & 1 & 1 & 1 & 1 & 1 & 1 & 1 & 1 & 0 & 0 & 0 & 0 & 1 & 1 & 1 & 1 & 0 &
0 & 0 & 0
\\
0 & 1 & 0 & 1 & 1 & 0 & 1 & 0 & 0 & 1 & 0 & 1 & 1 & 0 & 1 & 0 & 1 & 0 & 1 & 0 &
0 & 1 & 0 & 1 & 1 & 0 & 1 & 0 & 0 & 1 & 0 & 1 & 0 & 1 & 0 & 1 & 1 & 0 & 1 & 0 & 0
& 1 & 0 & 1 & 1 & 0 & 1 & 0 & 1 & 0 & 1 & 0 & 0 & 1 & 0 & 1 & 1 & 0 & 1 & 0 & 0 &
1 & 0 & 1
\\
\end{array} %
\rright).}}
\]

The constants for calculating wordlengths are
$1,1,1,2,2,2,2,2,2,2,2,2,3,\break 3,3,3,3,3,3,3,3,3,4,4,4,4,4,4,4,4,4,5,5,5,6$
and $\delta$ for seven a-equations are $0,0,0,1,1,1,0$. Then
we apply the $K$ and $A$ matrix to the design properties of a general
$(1/64)${th}-fraction QC design $D$ with an even number of factors.
Assume $D$ is constructed from a generator matrix $G=(u,v,w,I_n)$.
Theorem~\ref{th3} presented below gives an account of words of all possible types.

\begin{theorem}\label{th3} 
With reference to the $2^{(2n+6)-6}$ QC design $D$, assuming $\sum_{i=1,3;j=1,3;k=0,2} f_{ijk}$, $\sum_{i=1,3;j=0,2;k=1,3} f_{ijk}$ and
$\sum_{i=0,2;j=1,3;k=1,3} f_{ijk}$ are all greater than $0$, the
following hold:
\begin{longlist}[(a)]
\item[(a)] There are $8/\rho_{100}^2$ words each with aliasing index
$\rho_{100}$; each $1/4$ of them have lengths $k_{100}+1$, $k_{120}+3$,
$k_{102}+3$ and $k_{122}+5$.
\item[(b)] There are $8/\rho_{010}^2$ words each with aliasing index
$\rho_{010}$; each $1/4$ of them have lengths $k_{010}+1$, $k_{210}+3$,
$k_{012}+3$ and $k_{212}+5$.
\item[(c)] There are $8/\rho_{001}^2$ words each with aliasing index
$\rho_{001}$; each $1/4$ of them have lengths $k_{001}+1$, $k_{201}+3$,
$k_{021}+3$ and $k_{221}+5$.
\item[(d)] There are $8/\rho_{110}^2$ words each with aliasing index
$\rho_{110}$; each $1/4$ of them have lengths $k_{110}+2$, $k_{130}+2$,
$k_{112}+4$ and $k_{132}+4$.
\item[(e)] There are $8/\rho_{101}^2$ words each with aliasing index
$\rho_{101}$; each $1/4$ of them have lengths $k_{101}+2$, $k_{103}+2$,
$k_{121}+4$ and $k_{123}+4$.
\item[(f)] There are $8/\rho_{011}^2$ words each with aliasing index
$\rho_{011}$; each $1/4$ of them have lengths $k_{011}+2$, $k_{013}+2$,
$k_{211}+4$ and $k_{213}+4$.
\item[(g)] There are $8/\rho_{111}^2$ words each with aliasing index
$\rho_{111}$; each $1/4$ of them have lengths $k_{111}+3$, $k_{113}+3$,
$k_{131}+3$ and $k_{133}+3$.
\item[(h)] There are $7$ words each with aliasing index $1$; they have
lengths $k_{200}+2$, $k_{020}+2$, $k_{002}+2$, $k_{220}+4$,
$k_{202}+4$, $k_{022}+4$ and $k_{222}+6$, respectively.
\end{longlist}
All $\rho_{ijk}$ are defined as $2^{-\lfloor(a_{ijk}+\delta)/2 \rfloor
}$, where $\delta=1$ for $\rho_{110}$, $\rho_{101}$ and $\rho_{011}$,
and $\delta=0$ otherwise.
\end{theorem}

The proof of Theorem~\ref{th3} can be done in a similar way as either the
matrix expansion method in the proof of Theorem 1 of \citet{PhoXu09}
or the trigonometric approach in the proof of Theorem 2 of \citet{Zhaetal11}
and omitted here. Theorem~\ref{th3}, in conjunction with equations of
$K$ and $A$, shows that the resolution and wordlength pattern of the
design $D$ depend on $u$, $v$ and $w$ only. The following example
illustrates the calculations of the generalized resolution and
generalized wordlength pattern of $D$.

\begin{example}\label{ex4} 
Given the generating matrix of a $256 \times14$ quaternary-code design $D,$
\[
G=(u,v,w,I_4)= %
\pmatrix{ 1 & 1 & 2 & 1 & 0 & 0 & 0
\vspace*{2pt}
\cr
1 & 2 & 1 & 0 & 1 & 0 & 0 \vspace*{2pt}
\cr
1 & 3 & 3 & 0 & 0 & 1
& 0 \vspace*{2pt}
\cr
2 & 1 & 3 & 0 & 0 & 0 & 1 }. %
\]
$D$ can be represented by a frequency vector $F=(\vec{0}_{22},1,\vec
{0}_{2},1,\vec{0}_{5},1,\vec{0}_{7},1,\vec{0}_{24})$, where $\vec{0}_n$
is a vector of $0$ with length $n$. So $K=CF=( 5, 5, 5, 6, 4, 4, 6, 4,
4, 4, 4,\break  6, 3, 3, 3, 3, 3, 3, 7, 3, 3, 3, 4, 2, 6, 2, 6, 4, 6, 2, 4, 5,
5, 5, 2)$ and $A=BF=( 3, 3, 3,  2,\break 2, 2, 1 )$. It leads to 35
wordlengths with lengths 6, 6, 6, 8, 6, 6, 8, 6, 6, 6, 6, 8, 6, 6, 6,
6, 6, 6, 10, 6, 6, 6, 8, 6, 10, 6, 10, 8, 10, 6, 8, 10, 10, 10, 8 and seven
aliasing indexes all equal to $1/2$. Theorem~\ref{th3} entails 224 partial
words each with aliasing index $1/2$; of these, 168 have length six and
56 have length ten. In addition, Theorem~\ref{th3} entails seven complete words
of length eight. Hence, in this case the QC design $D$, which is a
$2^{14-6}$ design, has resolution $6.5$ and wordlength pattern
$(0,0,0,0,0,42,0,7,0,14,0,0,0,0)$. Comparing to the regular design of
the same size, this QC design has a higher resolution (6.5 versus 5.0)
and it has better aberration ($A_5=0$ for QC design versus $A_5 \neq0$
for regular design). Therefore, this QC design is more favorable than
its corresponding regular design.
\end{example}

Instead of performing a complete enumeration, a periodic structure for
a class of good $(1/64)${th}-fraction QC designs with high resolution
is presented in the following theorem.

\begin{theorem}\label{th4} 
Given a $2^{(2n+6)-6}$ QC design $D_0$ defined by a frequency vector
$F_0$, assume $D_0$ satisfies the conditions in Theorem~\ref{th3} and it has
generalized resolution $R_0=r_0+1-\rho_0$. Then for $t \geq0$, a
$2^{(2n+126t+6)-6}$ QC design $D_t$ defined by $F_t=F_0+(0,\vec
{1}_{63})t$ has generalized resolution $R_t=r_t+1-\rho_t$, where
$r_t=r_0+64t$ and $\rho_t=\rho_0(2^{-16t})$ if $\rho_0<1$ and $\rho_t=1$ if $\rho_0=1$.
\end{theorem}

\begin{example}\label{ex5} 
Following Example~\ref{ex4}, let $F_0 =(\vec{0}_{22},1,\vec{0}_{2},1,\vec
{0}_{5},1,\vec{0}_{7},1,\vec{0}_{24})$ and the $256 \times14$ QC
design $D_0$ has generalized resolution $6.5$. Then Theorem~\ref{th4} suggests
that for $t=1$, a $2^{140-6}$ QC design $D_1$ defined by $F_t =(0,\vec
{1}_{21},2,\vec{1}_{2},2,\vec{1}_{5},2,\vec{1}_{7},2,\vec{1}_{24})$ has
$r_t=6+64(1)=70$ and $\rho_t=(1/2)\times\break (2^{-16(1)})=2^{-17}$, that is,
generalized resolution $70.9999924$.
\end{example}

\section{Summary}\label{sec5}

This work provides some theoretical understandings of the structure of
a general $(1/4)^p$th-fraction QC design. In Section~\ref{sec2} we show via the
Code Arithmetic approach how the $k$-equations and a-equations of a
general $(1/4)^p$th-fraction QC design are developed from those of
other $(1/4)^h$th-fraction QC designs, where $p>h$. Section~\ref{sec3} lists
four rules on the structure of $k$-equations and a-equations when some
entries of $\vec{w}$ are added and/or changed. In addition, Theorem~\ref{th1}
describes the general structure of $k$-equations when all entries are odd
and Theorem~\ref{th2} suggests which $k$-equations are sufficient to be
considered so that the design properties can be determined. In Section
\ref{sec4} these rules and theorems are applied to determine the properties of
$(1/64)${th}-fraction QC designs and the periodic structure regarding
resolution is derived.

\section{Proofs}\label{sec6}

\subsection{\texorpdfstring{Proof of Corollaries \protect\ref{co1} and \protect\ref{co2}}
{Proof of Corollaries 1 and 2}}

We prove Corollary~\ref{co1} via induction. It is trivial for $p=1$, because it
leads to $k_{10}$ and $k_{01}$ for $l=0,1$. Assume $p=z$ is true, that
is, $k_{(\vec{0}_l,1,\vec{0}_{z-l})} = \sum_{\vec{i} \in C(z)} (f_{(\vec
{i}_{l},1,\vec{i}_{z-l})}+f_{(\vec{i}_{l},3,\vec{i}_{z-l})}+2 f_{(\vec
{i}_{l},2,\vec{i}_{z-l})})$. For $p=z+1$, we rewrite $\vec{w}$ as $(\vec
{0}_l,1,0,\vec{0}_{z-l})$, that is, insert a $0$ in the $(l+2)$th entry
of $\vec{w}$. Applying Rule~\ref{ru1}, we have $k_{(\vec{0}_l,1,0,\vec
{0}_{z-l})}=\sum_{s=0}^3 \sum_{\vec{i} \in C(z)} (f_{(\vec
{i}_{l},1,s,\vec{i}_{z-l})}+f_{(\vec{i}_{l},3,s,\vec{i}_{z-l})}+2
f_{(\vec{i}_{l},2,s,\vec{i}_{z-l})})$. Notice that $(\vec{i}_{l},s,\vec
{i}_{z-l})$ represents the $i${th} row of $C(z+1)$ for $s=0,1,2,3$,
and the above equation becomes $k_{(\vec{0}_l,1,\vec{0}_{(z+1)-l})}=
\sum_{s=0}^3 \sum_{\vec{i} \in C(z+1)} (f_{(\vec{i}_{l},1,\vec
{i}_{(z+1)-l})}+f_{(\vec{i}_{l},3,\vec{i}_{(z+1)-l})}+2 f_{(\vec
{i}_{l},2,\vec{i}_{(z+1)-l})})$. This completes the proof of Corollary
\ref{co1}. The proof of Corollary~\ref{co2} follows the same induction except the
formula is different.

\subsection{\texorpdfstring{Proof of Theorem \protect\ref{th1}}{Proof of Theorem 1}}

We prove Theorem~\ref{th1} via induction. The cases of $p=1$ and $p=2$ are true
from the results of \citet{PhoXu09} and \citet{Zhaetal11}. Assume
it is true for $p=z$ is true, that is, for $k_{\vec{w}} = 1 \sum_{\vec
{i} \in C_1(z)} f_{\vec{i}} + 2 \sum_{\vec{i} \in C_2(z)} f_{\vec{i}}$,
the sum of entries of all $\vec{i}$ in $C_1(z)$ are odd and the sum of
entries of all $\vec{i}$ in $C_2(z)$ are even. Consider $p=z+1$. We
start from rewriting $k_{\vec{1}_z} = 0(\sum_{\vec{i} \in C_0(z)}
f_{(\vec{i}_{z-1},i_{z})})+1(\sum_{\vec{i} \in C_1(z)} f_{(\vec
{i}_{z-1},i_{z})})+2(\sum_{\vec{i} \in C_2(z)} f_{(\vec
{i}_{z-1},i_{z})})$. The application of Rule~\ref{ru2} suggests that\break $
k_{(1,\vec{3}_z)} = \sum_{s=0}^3 0(\sum_{\vec{i} \in C_0(z)} f_{(s,\vec
{i}_{z-1},i_{z}+s)})+1(\sum_{\vec{i} \in C_1(z)} f_{(s,\vec
{i}_{z-1},i_{z}+s)})+\break2(\sum_{\vec{i} \in C_2(z)} f_{(s,\vec
{i}_{z-1},i_{z}+s)})$. Notice that if the sum of entries of $(\vec
{i}_{z-1})$ plus $i_z$ is odd, then $s$ plus the sum of entries of $\sum (\vec{i}_{z-1})$
plus $(i_z+s)$ is still odd for $s=0,1,2,3$. It is
also true for the even case.

Applying Rule~\ref{ru3}, $k_{(\vec{w}_l,s_2,\vec{w}_{p-l})}=k_{(\vec
{w}_l,s_1,\vec{w}_{p-l})} \oplus k_{(\vec{0}_l,2,\vec{0}_{p-l})}$,
where $s_2=(s_1+2) \operatorname{mod} 4$. $L_w(1+2)=L_w(3+2)=1$ implies that the
frequencies with an odd sum of entries of $\vec{i}$ have odd
coefficients. Similarly, $L_w(0+2)=2$ and $L_w(2+2)=0$ imply that the
frequencies with an even sum of entries $\vec{i}$ have even
coefficients. Therefore, by repeatedly applying Rule~\ref{ru3} to change all
entries of $3$ into $1$ in $\vec{w}$, we can express $k_{\vec{1}_{z+1}}
= 0\sum_{\vec{i} \in C_0(z+1)} f_{\vec{i}}+1\sum_{\vec{i} \in C_1(z+1)}
f_{\vec{i}}+2\sum_{\vec{i} \in C_2(z+1)} f_{\vec{i}}$. This completes
the proof.

\subsection{\texorpdfstring{Proof of Theorem \protect\ref{th2}}{Proof of Theorem 2}}

Consider a general $(1/4)^p$th-fraction QC design $D$. There are $4^p$
different combinations of $\vec{w}$ with entries in $Z_4 \in\{0,1,2,3\}
$. Among these $\vec{w}$, there are $2^p$ of them where their entries
are all even. Then it is obvious that $k_{\vec{0}}$ is obviously
irrelevant to any properties of $D$ because this $k$-equation does not
include any columns from $V$ and the columns from $I_n$ are complete.
This leads to the first group of $\vec{w}$ with a total of $2^p-1$
possible combinations.

Eliminating the choice with all even entries, there are $4^p-2^p$
different $\vec{w}$ that consist of at least one odd entry. If we focus
on the first odd entry of $\vec{w}$, half of these $\vec{w}$ start with
$1$ and another half start with $3$. Notice that $k_{\vec{w}}$ and
$k_{\vec{w'}}$ are equivalent if all $1$ entries in $\vec{w}$ become
$3$ entries in $\vec{w'}$ and vice versa. It is proved as follows.

Using the expression in Theorem~\ref{th1}, without loss of generality, $k_{\vec
{w}} = \break 1 \sum_{\vec{i} \in C_1(p)} f_{\vec{i}} + 2 \sum_{\vec{i} \in
C_2(p)} f_{\vec{i}}$. A repeated use of Rule~\ref{ru3} on every odd entry of
$\vec{w}$ in $k_{\vec{w}}$ leads to $k_{\vec{w'}} = k_{\vec{w}} \oplus
k_{\vec{w}_2}$, where the entries of $\vec{w}_2$ are $2$ if the
corresponding entry of $\vec{w}$ is odd, and $0$ otherwise. We can
express $k_{\vec{w}_2}$ easily by the CA operation on the expressions
of Corollary~\ref{co2}\vspace*{1pt} and it results in $k_{\vec{w}_2} = 2 \sum_{\vec{i} \in
C_1(p)} f_{\vec{i}} + 0 \sum_{\vec{i} \in C_2(p)} f_{\vec{i}}$. Then
$k_{\vec{w'}}$ can be expressed in the same way as $k_{\vec{w}}$ due to
the Lee weight $L_w(3)=1$.

Therefore, for all $\vec{w}$ that consist of odd entries, it is
sufficient and necessary to consider the $k$-equations that the first odd
entry of $\vec{w}$ is $1$, and there are $(4^p-2^p)/2$ or
$2^{2p-1}-2^{p-1}\vec{w}$ in total.

\subsection{\texorpdfstring{Proof of Theorem \protect\ref{th4}}{Proof of Theorem 4}}
About the periodicities of $r_t$, we start from\vspace*{1pt} the original $k$-matrix
$K_0=CF_0$. If $F_t=F_0+(0,\vec{1}_{63})t$, then $K_t=CF_t=C(F_0+(0,\vec{1}_{63})t)=K_0+C(0,\vec{1}_{63})t$. Since the second
term results in a vector of length~$35$ and all entries are $64t$, and
the constants for calculating wordlengths are invariant to $t$, $r_t = r_0+64t$.

About the periodicities of $\rho_t$, we start from the original
a-matrix $A_0=BF_0$. Similar to the $k$-matrix, $A_t=A_0+B(0,\vec
{1}_{63})t$. Since the second term results in a vector of length $7$
and all entries are $32t$, and the constants for calculating aliasing
indexes are fixed at $(0,0,0,1,1,1,0)$, $\rho_t=2^{-\lfloor(a_t+\delta
)/2 \rfloor}=2^{-\lfloor(a_0+32t+\delta)/2 \rfloor}=2^{-\lfloor
(a_0+\delta)/2 \rfloor}2^{-32t/2}=\rho_0(2^{-16t})$.

\section*{Acknowledgments}
The author
would like to thank the Associate Editor, two referees and Professor
Hongquan Xu for their valuable suggestions and comments to this paper.


%


%
%

\printaddresses

\end{document}